\documentstyle[12pt]{article}
\input mssymb.tex
\newtheorem{thm}{Theorem}[section]
\newtheorem{prop}[thm]{Proposition}
\newtheorem{define}[thm]{Definition}
\newtheorem{cor}[thm]{Corollary}
\newtheorem{lem}[thm]{Lemma}

\begin{document}

\parskip0.2cm

\title{Linearization of Nambu structures}

\author{Jean-Paul Dufour and Nguyen Tien Zung \\
D\'epartement de Math\'ematiques, Universit\'e Montpellier II\\
}

\date{7 July 1997}

\newpage
\maketitle

\begin{abstract} {\normalsize
Nambu structures are a generalization of Poisson structures in Hamiltonian 
dynamics, and it has been shown recently by several authors that, outside 
singular points, these structures are locally an exterior product of commuting vector 
fields. Nambu structures also give rise to co-Nambu differential forms, which 
are a natural generalization of integrable 1-forms to higher orders. This work 
is devoted to the study of Nambu tensors and co-Nambu forms near singular 
points. In particular, we give a classification of linear Nambu structures 
(also called finite-dimensional Nambu-Lie algebras), and a linearization of 
Nambu tensors and co-Nambu forms, under the nondegeneracy condition. 


}

\end{abstract}



\section{Introduction}

Let $V$ be an $n$-dimensional smooth manifold and $C = C^{\infty}(V)$
the space of smooth functions on $V$. A Nambu structure of order $q$ on V
is a multi-linear anti-symmetric application $\Pi$ from	the direct product
of $q$ samples of $C$ to $C$, and denoted by the bracket $\{ \}$:
$$
\Pi: C \times C \times \ldots \times C \longrightarrow C, \;
(f_1, f_2, \ldots, f_q) \mapsto \{f_1, f_2, \ldots, f_q \} 
$$
which satisfies the following two conditions \\

i) Leibnitz condition:
\begin{equation} 
\label{eq:Leibnitz}
\Pi_{f_1,...,f_{q-1}} (fg) = f\Pi_{f_1,...,f_{q-1}} (g) + g \Pi_{f_1,...,f_{q-1}}(f)
\end{equation}

ii) Jacobi condition:
\begin{equation} 
\label{eq:Jacobi}
\Pi_{f_1,...,f_{q-1}} (\{g_1,...,g_q\}) = \sum^{q}_{i=1} 
\{ g_1,..., g_{i-1}, \Pi_{f_1,...,f_{q-1}} (g_i),..., g_q \}
\end{equation}
for any $f_1,...,f_{q-1},g_1,...,g_{q-1},f,g \in C$, where $\Pi_{f_1,...,f_{q-1}}$
denotes the contraction of $\Pi$ by $f_1,...,f_{q-1}$:
$\Pi_{f_1,...,f_{q-1}}(f) := \{ f_1,...,f_{q-1},f \}$

The Leibnitz condition (together with the antisymmetricity of $\Pi$)
means that $\Pi$ is given by an (anti-symmetric) $q$-vector field on $V$,
which we will also denote by $\Pi$. When $q=1$ the Jacobi condition is empty
and we simply have a vector field on $V$. When $q=2$ the Jacobi condition
is the usual condition for a 2-vector field to be a Poisson structure in
Hamiltonian dynamics. Thus, Nambu structures, which are also called
Nambu-Poisson structures, are a kind of generalization of Poisson structures
when the order $q$ is different from 2. 
They were introduced by Nambu \cite{Nambu}
in an attempt to generalize Hamiltonian mechanics.

Given a Nambu structure of order $q$ and a $(q-1)$-tuple of functions
$(f_1,...,f_{q-1})$ on $V$, one can associate to it a {\it Hamiltonian}
vector field, which is the vector field corresponding to the derivation
$\Pi_{f_1,...,f_{q-1}}: C \to C$. The Jacobi condition 
means that this Hamiltonian vector field preserves the Nambu structure,
like in usual Hamiltonian dynamics. From the definition it is evident that
the contraction $\Pi_{f_1,...,f_{q-r}}$ of a Nambu structure $\Pi$
of order $q$ with arbitrary $q-r$ smooth functions $f_1,...,f_r$
($0 < r < q$),$\Pi_{f_1,...,f_{q-r}}(g_1,...,g_{r}) :=          
\Pi_{f_1,...,f_{q-r}}(g_1,..., g_{r})$
is again a Nambu structure of order $r$. In particular, when $q \geq 3$
and $r =2$, we get an infinite family of Poisson structures.

Nambu structures were studied by many people in recent years, and one
can imagine various algebraic structures associated to them (\cite{Gautheron,Takhtajan}).
The most significant result obtained, which is in fact also quite simple to prove,
is the following local normal form theorem, which was proved by Gautheron \cite{Gautheron}
and independently by Nakanishi \cite{Nakanishi}, Alekseevsky and Guha \cite{AG}. 
Hereafter by a {\it Nambu tensor} of order $q$ we will mean an $q$-vector
field associated to a Nambu structure.

{\bf Theorem (Gautheron et al.)}.
{\it
Let $\Pi$ be a Nambu tensor of order $q \geq 3$ on an $n$-dimensional manifold
$V$, and $O \in V$ a point in which $\Pi(O) \neq 0$.
Then in a small neighborhood of $O$ one can find a local system
of coordinates $(x_1,...,x_n)$ such that
$\Pi = \partial / \partial x_1 \wedge ... \wedge \partial / \partial x_q$
in this neighborhood.
}

The above theorem is a kind of Darboux theorem for Nambu structures.
It also shows a big difference between Nambu structures of order $ \geq 3$
and Poisson structures: the former ones are decomposable at non-zero points
while the later ones are not in general.

The above theorem prompts us to study singularities of Nambu structures.
The first obvious thing that we observe here is that each Nambu structure gives rise
to an associated singular foliation (in the sense of Stefan-Sussmann), whose
distribution is spanned by the Hamiltonian vector fields $\Pi_{f_1,...,f_{q-1}}$.
When $q \geq 3$ the leaves of this singular foliation is of dimension either
0 or $q$, while in case  of Poisson structures ($q =2$) they may have any
even dimension (see e.g. \cite{Vaisman,Weinstein} for the case of Poisson structures).
These singular foliations give a geometric picture about the Nambu structures 
themselves.

By a {\it singularity of a Nambu structure $\Pi$}, or a {\it Nambu singularity}  
we mean a small neighborhood of a point $O$ at which	$\Pi(O) = 0$.
When $\Pi(O) = 0$ at some point $O$ , then its linearization at $O$ is well-defined
and gives us a linear Nambu structure. Thus the study of linear Nambu structures
is a natural first step in the study of singularities of general Nambu
structures. We have the following result (cf. Corollary \ref{cor:LinearTensor})

\begin{thm}                                                                  
Every linear Nambu tensor $\Pi$ of order $q = n-p \geq 3$ on an             
$n$-dimensional linear space $V$ belongs to one of the following two types:  
                                                                             
\underline{Type 1}:                                                          
$                                                                            
\Pi = \sum_{j=1}^{r+1} \pm x_j \partial / \partial x_1 \wedge ... \wedge     
\partial / \partial x_{j-1}	\wedge                                         
\partial / \partial x_{j+1}	\wedge ... \wedge                              
\partial / \partial x_{q+1} +   \\                                             
\sum_{j=1}^{s} \pm x_{q+1+j} \partial / \partial x_{1} \wedge ... \wedge     
\partial / \partial x_{r+j} \wedge                                           
\partial / \partial x_{r+j+2} \wedge                                         
\partial / \partial x_{q+1}                                                  
$                                                                            
(with $-1 \leq r \leq q, 0 \leq s \leq min(p-1, q-r)$).                      
                                                                             
\underline{Type 2}:                                                          
$                                                                            
\Pi =                                                                        
\partial / \partial x_{1} \wedge ... \wedge                                  
\partial / \partial x_{q-1} \wedge                                           
(\sum_{i,j=q}^{n} b^i_j x_i \partial / \partial x_j)                        
$                                                                           
\end{thm}          

We will call a Nambu singularity {\it of Type 1} if its linear part           
is of Type 1, and {\it of Type 2} in the other case. The singularities of     
Type 1 and Type 2 are very different geometrically, their corresponding       
foliations look very different, though they are in some	natural             
sense {\it dual} to each other (cf. Section \ref{section:Linear}). We have the 
following result about the linearization of Nambu tensors near singular points
(see Theorem \ref{thm:FormType1},
Theorem \ref{thm:LinearizationType1} and Theorem \ref{thm:LinearizationType2}
for the precise formulations):

\begin{thm}
Nondegenerate singularities of Type 1 of Nambu tensors of order $q \geq 3$
are formally linearizable. They are, up to multiplication by a function, 
$C^{\infty}$-linearizable if they are elliptic, 
and $C^{\omega}$-linearizable in the analytic (real or complex) case. 
Nondegenerate singularities of Type 2 of Nambu tensors of order $q \geq 3$ are 
$C^{\infty}$-linearizable under some nonresonance condition, and
analytically linearizable in the analytic case under some Diophantine condition.
\end{thm}

For non-elliptic singularities of type 1 of 
class $C^{\infty}$, we have (see Section 
\ref{section:Type1}): In the case of signature $q-3$
they are not continuously linearizable in general. 
If the signature is different from $q-3$ 
then they are conjectured to be $C^{\infty}$-linearizable. What we know is that
in this case their associated singular foliations are homeomorphic to the ones 
given by the linear Nambu structures.

An important object which arises in the study of Nambu tensors are the 
so-called {\it co-Nambu forms}, which are 
obtained by the contraction of Nambu tensors with volume forms. For them we 
have some results analogous to the above theorem, 
which complement the ones obtained by Medeiros \cite{Medeiros}, 
and are similar to some
results obtained before by Kupka \cite{Kupka}, Reeb \cite{WR}, Moussu 
\cite{Moussu,Moussu2} and others for integrable 1-forms. 
Thus one can think of co-Nambu forms as 
integrable differential forms of higher orders. In fact, they are called 
{\it integrable p-forms} in \cite{Medeiros}. In particular, we suspect that 
many results obtained by various authors
for degenerate singularities of integrable 1-forms can be also generalized
to the case of co-Nambu forms. 

The rest of this paper is organized as follows: In Section 
\ref{section:Preliminaries} we give some preliminary results concerning 
Nambu structures, most notably about co-Nambu forms. In Section 
\ref{section:Linear} we give a classification of linear Nambu structures, where
we show that they can be divided in two types. In Section 
\ref{section:Decomposition} we prove a theorem about the decomposition of 
Nambu structures near nondegenerate singularities. 
Section \ref{section:Type1} and Section \ref{section:Type2} 
contain our main results concerning the linearization problem.

{\bf Acknowledgements}. We thank N. Nakanishi whose visit to Montpellier in
May 1997 attracted our attention to this problem. 

\section{Preliminaries}
\label{section:Preliminaries}

Let $\Omega$ be a volume form on an $n$-dimensional manifold  $V$, and
$\Pi$ an $q$-vector field on $V$, with 
$n > q > 2$. Put $p = n -q $ and denote by $\omega$ the $p$-form obtained 
by contracting $\Pi$ and $\Omega$:
$$ \omega = i_{\Pi} \Omega $$
Then the condition for $\Pi$ to be a Nambu tensor can be rewritten
in terms of $\omega$:

\begin{prop}
\label{prop:omega}
With the above notations, $\Pi$ is Nambu if and only if
$\omega$ satisfies the following two conditions: \\
\begin{equation}
\label{eq:Decomposable}
 i_A \omega \wedge \omega = 0 
\end{equation}
\begin{equation}
\label{eq:Integrable}
 i_A \omega \wedge d\omega = 0
\end{equation}
for any $(p-1)$-vector $A$.
\end{prop}

In case $p = 1$ the above conditions simply mean that 
$d\omega \wedge \omega = 0$, i.e. $\omega$ is an integrable 1-form.

The proof of the above proposition is based on the following two lemmas,
which follow directly from the Leibnitz and Jacobi conditions 
(\ref{eq:Leibnitz}), (\ref{eq:Jacobi}) and the normal form theorem of Gautheron et 
al.

\begin{lem} 
\label{lem:1}
$\Pi$ is a Nambu tensor if and only if
it is so on the open set $U = \{ x \in V, \Pi(x) \neq 0\}$
of points where it does not vanish.
\end{lem}

\begin{lem}
\label{lem:2}
Suppose that $ q \geq 3$. Then a $q$-vector field $\Pi$ is Nambu
if and only if in a neighborhood of each point $O$ where $\Pi(O) \neq 0$, 
we can find a local system of coordinates
$(x_1,...,x_n)$ in which $\Pi$ can be written as 
$\Pi = \partial/\partial x_1 \wedge \ldots \wedge \partial / \partial x_q $.
\end{lem}

{\it Proof of Proposition \ref{prop:omega}}. Let $\Pi$ be a Nambu
$q$-tensor with $q \geq 3$. In a neighborhood of a point $O$
such that $\Pi(O) \neq 0$ we have 
$\Pi = \partial/\partial x_1 \wedge \ldots \wedge \partial / \partial x_q $
in some system of coordinates, according to the theorem of Gautheron et al..
Since $\Omega = f dx_1 \wedge ... \wedge dx_n$ (with some non-zero function
$f$), we have
$$ \omega = \pm f dx_{q+1} \wedge ... \wedge dx_n $$
and 
$$ d\omega = \pm df \wedge dx_{q+1} \wedge ... \wedge dx_n$$

From here it is easy to verify that the Equation (\ref{eq:Decomposable})
and Equation (\ref{eq:Integrable}) 
are satisfied for any $(p-1)$-vector $A$, where
$p = n - q$. (At least they are satisfied at any non-zero point of $\Pi$, 
but then at any point, since zero points of $\Pi$ are also zero points
of $\omega$).

Conversly, let $\Pi$ be a $q$-vector such that 
$\omega = i_{\Pi}\Omega$ satisfies the equations (\ref{eq:Decomposable})
and (\ref{eq:Integrable}). Fix a point $O \in V$ such that
$\Pi(O) \neq 0$ (hence $\omega(O) \neq 0$). Then
Equation (\ref{eq:Decomposable}) implies that $\omega$ is decomposable
in a neighborhood of $O$:
$$ \omega = \alpha_1 \wedge ... \wedge \alpha_p$$
where $\alpha_i$ are independent 1-forms. One can find
$(p-1)$-vectors $A_1,..., A_p$ such that
$$i_{A_j} \omega = \alpha_j, j = \overline{1,p} $$
in some neighborhood of $O$. Hereafter $\overline{1,p}$ means $1,2,...,p$.
Then Equation (\ref{eq:Integrable})	gives
$ \alpha_j \wedge d\omega = 0 $ for  $ j=\overline{1,p} $
in this neighborhood. But
$$ d\omega = \sum_{j=1}^p \alpha_1 \wedge ... \wedge \alpha_{j-1} \wedge d \alpha_j
\wedge \alpha_{j+1} \wedge ... \wedge \alpha_p $$
Thus we have
$$ d\alpha_j \wedge \alpha_1 \wedge ...\wedge \alpha_p = 0 \; \; {\rm for}\; j = \overline{1,p}$$
In other words, $\alpha_j$ satisfy the conditions of Frobenius theorem 
(see e.g. \cite{BCGGG}), which says that in this case there exists a local system of 
coordinates $(x_1,...,x_n)$ such that 
$\alpha_j \wedge dx_{q+1} \wedge ... \wedge dx_n = 0$   ${\rm for}\; j = \overline{1,p}$.
It follows that $\omega = f dx_1 \wedge ... \wedge dx_p$ for some non-zero
function $f$, 
and $\Pi = g \partial / \partial x_1 \wedge ... \wedge \partial / \partial x_q$
for some non-zero function $g$. Replacing $x_1$ by 
$x_1' = \int_{t=0}^{x_1} \frac{dt}{g(t,x_2,...,x_n)}$, we have
$\Pi = \partial / \partial x_1 \wedge ... \wedge \partial / \partial x_q$.
Applying Lemma \ref{lem:2}, we obtain that $\Pi$ is a Nambu structure 
$\Box$

A simple corollary of Proposition \ref{prop:omega} is that if $\Pi$ is a Nambu
tensor of order $q \geq 3$ and if $f$ is a smooth function, then $f\Pi$ is again
a Nambu structure. 

\begin{define}
\label{define:coNambu}
A differential $p$-form $\omega$ which satisfies the equations 
(\ref{eq:Decomposable})
and (\ref{eq:Integrable}) in Proposition \ref{prop:omega} will be called a {\bf
co-Nambu form} (of order $p$ and co-order $q$).
\end{define}

We have a bijection $\Pi \leftrightarrow \omega$ between Nambu tensors and
co-Nambu forms (if $V$ is orientable). Of course, this bijection depends on 
the choice of a volume form on $V$, so it is not unique, but unique up to
multiplication by a non-zero function. Thus the study of singularities of
$\Pi$ and that of $\omega$ are almost the same.

As a principle, when a structure vanishes at some point, then its linearization
is well-defined, and if its linearization also vanishes, then its quadratization
is well-defined, etc. It is also true for Nambu and co-Nambu structures.
Let $O \in V$ be a point such that $\Pi(O) = 0$, and $(x_1,...,x_n)$ a local system
of coordinates in a neighborhood of $O$. Then we have a Taylor expansion of $\Pi$
at $O$: 
$$ \Pi = \Pi^{(1)} +  \Pi^{(2)} + \Pi^{(3)} + ... $$
where 
$$  \Pi^{(i)} = \sum_{j_1 \leq ... \leq j_q} P^{(i)}_{j_1 ... j_q}
\partial / \partial x_{j_1} \wedge ... \wedge \partial / \partial x_{j_q}$$
with $P^{(i)}_{j_1 ... j_q}$ being polynomials of order $i$ in $x_1,...,x_n$.
It is easy to see from the definition that $\Pi^{(1)}$ is well-defined, and is also a Nambu 
structure. It is called the {\it linear part of $\Pi$}.

Similarly (by putting $\Omega = dx_1 \wedge ... \wedge dx_n$) we have
$
\omega =  \omega^{(1)} +  \omega^{(2)} + \omega^{(3)} + ... 
$,
with
$
\omega^{(k)} = i_{\Pi^{(k)}} \Omega =
\sum_{j_1 \leq ... \leq j_q} \pm P^{(k)}_{j_1 ... j_q} dx_1 \wedge ... \wedge
\hat{dx_{j_1}} \wedge ... \wedge \hat{d x_{j_q}} \wedge ... \wedge dx_n 
$.
In particular, the linear part $\omega^{(1)}$ of $\omega$ is well-defined by $\omega$
and is also a co-Nambu form. Note that 	$\omega^{(1)}$ is uniquely determined by 
$\Pi^{(1)}$, up to multiplication by a constant.

For co-Nambu 1-forms, Proposition \ref{prop:omega} shows that they are nothing but 
integrable 1-forms. (This has been known to be true also for Poisson structures on 
3-manifolds, cf. \cite{Dufour}). The singularities of integrable  1-forms have been 
extensively studied (see e.g. \cite{Kupka,Malgrange,Medeiros,Moussu,WR}). 
In particular, there is the following so-called Kupka's phenomenon 
(see \cite{Kupka,Medeiros}): If $O$ is a zero point of an integrable 
1-forms $\omega$ and $d\omega(O) \neq 0$, then locally $\omega$ is a pull-back 
of an 1-form on a plane. In \cite{Medeiros} a similar result is also proved for
co-Nambu forms of higher orders.

\section{Linear Nambu structures}
\label{section:Linear}

\begin{thm}
\label{thm:Linear}
If $\omega$ is a linear co-Nambu $p$-form of co-order $q = n-p \geq 3$ on a linear space
$V$ then there exist linear coordinates $(x_1,...,x_n)$ such that $\omega$ belongs
to one of the following two types:

\underline{Type 1}: $ \omega = dx_1 \wedge ... \wedge dx_{p-1} \wedge \alpha$,
where $\alpha$ is an exact 1-form of the type
$
\alpha = d [ \sum_{j=p}^{p+r} \pm x_j^2 / 2 + \sum_{i=1}^s x_i x_{p+r+i}], 
$
with $ - 1 \leq r \leq q = n-p, 0 \leq s \leq q - r$.

\underline{Type 2}: $ \omega = \sum_{i=1}^{p+1} a_i dx_1 \wedge ... \wedge dx_{i-1} \wedge 
dx_{i+1} \wedge ... \wedge dx_{p+1}$
with $a_i = \sum_{j=1}^{p+1} a_i^j x_j$, where $a^i_j$ are constant.
The matrix $(a_i^j)$ can be chosen to be in Jordan form.
\end{thm}

{\it Proof.} Put $\omega = \sum_{j=1}^{n} x_j \omega_j$ where $\omega_j$ are constant 
$p$-forms. Then $\omega = \omega_j$ at points $(x_1 = 0,..., x_j = \epsilon,..., x_n = 0 )$.
At any point $\omega$ is either decomposable (i.e. a wedge product of covectors) or zero,
so does $\omega_j$ since it is constant. Denote by $E_j$ the span of $\omega_j$, i.e.
$$
\begin{array}{c}
E_j = Span (\omega_j) 
       \stackrel{def}{=} Span \{i_A \omega_j ,\; A \; is \; a \; 
       (p-1)-covector\} \\
	 = Annulator\{x \in V, \, i_x \omega_j = 0 \} \subset V
\end{array}
$$
Then $\dim E_j = p$ if $\omega_j \neq 0$, because of decomposability. We have:

\begin{lem}
\label{lem:p-1}
If $\omega_i \neq 0$ and $\omega_j \neq 0$ for some indices $i$ and $j$, then
$\dim (E_i \cap E_j) \geq p-1$.
\end{lem}

{\it Proof of Lemma \ref{lem:p-1}}. Putting $x_k = 0$ for every 
$k\neq i,j$, we obtain that
$x_i\omega_i + x_j \omega_j = \omega$ is decomposable or null for any $x_i,x_j$. 
In particular, $\omega_i + \omega_j$ is decomposable. If $\dim (E_i \cap E_j) = d < p$
then there is a basis \\ $(e_1,...,e_d,f_1,...,f_{p-d},g_1,...,g_{p-d})$ of $E_i + E_j$
such that $\omega_i = e_1 \wedge ... \wedge e_d \wedge f_1 \wedge ... \wedge 
f_{p-d}$, $\omega_j = e_1 \wedge ... \wedge e_d \wedge g_1 \wedge ... \wedge 
g_{p-d}$ and 
$$
\omega_i + \omega_j = e_1 \wedge ... \wedge e_d \wedge 
[f_1 \wedge ... \wedge f_{p-d} + \wedge g_1 \wedge ... \wedge g_{p-d}]
$$
It follows easily that if $p-d \geq 2$ then $Span(\omega_i + \omega_j) = E_i + E_j$,
$\dim Span(\omega_i + \omega_j) > p $ and $\omega_i + \omega_j$ is not decomposable.
$\Box$

Return now to Theorem \ref{thm:Linear}.  We can assume that $E_1,...,E_h \neq 0$ and
$E_{h+1},..., E_n = 0$ for some number $h$. Put $E = E_1 \cap E_2 \cap ... \cap E_h$.
Then there are two alternative cases: $\dim E \geq p-1$ and $\dim E < p-1$.

\underline{Case 1}. $\dim E \geq p-1$. Then denoting by $(x_1, ..., x_{p-1})$ a set of $p-1$
linearly independent covectors contained in $E$, and which are considered as linear 
functions on $V$, we have 
$$\omega_i = dx_1 \wedge dx_2 \wedge ... \wedge dx_{p-1} \wedge \alpha_i, \, i=\overline{1,h}$$
for some constant 1-forms $\alpha_i$, and hence
\begin{equation}
\label{eq:TempoType1}
\omega = dx_1 \wedge dx_2 \wedge ... \wedge dx_{p-1} \wedge \alpha
\end{equation}
where $\alpha = \sum x_i \alpha_i$ is a linear 1-form.

\underline{Case 2}. In this case, without loss of generality, we can assume that
$\dim (E_1 \cap E_2 \cap E_3) < p-1$. Then Lemma \ref{lem:p-1} implies that
$\dim (E_1 \cap E_2 \cap E_3) = p-2$. For an arbitrary index $i$, $3 < i \leq h$, 
put $F_1 = E_1 \cap E_i, F_2 = E_2 \cap E_i, F_3 = E_3 \cap E_i$. Recall that
$\dim F_1, \dim F_2, \dim F_3 \geq p-1$ according to Lemma \ref{lem:p-1}, but
$\dim (F_1 \cap F_2 \cap F_3) = \dim(E_1 \cap E_2 \cap E_3 \cap E_i) < p-1$, hence 
we cannot have $F_1 = F_2 = F_3$. Thus we can assume that
$F_1 \neq F_2$. Then either $F_1$ and $F_2$ are two different hyperplanes in $E_i$, 
or one of them coincides with $E_i$. In any case
we have $E_i = F_1 + F_2 \subset E_1 + E_2 + E_3$. It follows that
$$\sum_1^n E_i = \sum_1^h E_i = E_1 + E_2 + E_3$$
On the other hand, we have $\dim (E_1 + E_2 + E_3) = \dim E_1 + \dim E_2 + \dim E_3
- \dim(E_1 \cap E_2) - \dim (E_1 \cap E_3) - \dim (E_2 \cap E_3) + \dim (E_1 + E_2 + E_3)
= 3p - 3 (p-1) + (p-2) = p+1$. Thus
$$ \dim (E_1 + E_2 + ... + E_n) = p+1$$
It follows that there is a system of linear coordinates $(x_1,...,x_n)$ on $V$ such that
$(x_1,...,x_{p+1})$ span $E_1 + ... + E_n$ and therefore
$$ \omega_i = \sum_{j=1}^{p+1} \gamma_i^j dx_1 \wedge ... \wedge dx_{j-1} 
\wedge dx_{j+1} \wedge ... \wedge x_{p+1}$$
Hence we have
\begin{equation}
\label{eq:TempoType2}
\omega = \sum x_i \omega_i = \sum_{j=1}^{p+1} a_j dx_1 \wedge ... \wedge dx_{j-1}
\wedge dx_{j+1} \wedge ... \wedge x_{p+1}
\end{equation}
where $a_j $ are linear functions on $V$.

To finish the proof of Theorem \ref{thm:Linear}, we still need to normalize further
the obtained forms (\ref{eq:TempoType1}) and (\ref{eq:TempoType2}). 

Return now to Case 1 and suppose that 
$\omega = dx_1 \wedge ... \wedge dx_{p-1} \wedge \alpha$
where $\alpha = \sum \alpha_j dx_j$ with $\alpha_j$ being linear functions. We can put
$\alpha_j = 0$ for $j = \overline{1, p-1}$ since it will not affect $\omega$. 
Then we have $\alpha = \sum_{j \geq p, i = \overline{1,n}} \alpha^i_j x_i dx_j$. 
Equation (\ref{eq:Integrable}) implies
that $\alpha \wedge dx_1 \wedge ...\wedge dx_{p-1} \wedge d\alpha = 0$. If we consider 
$(x_1,...,x_{p-1})$ as parameters and denote by $d'$ 
the exterior derivation with respect to the variables $(x_p,...,x_n)$, 
then the last equation means $\alpha \wedge d'\alpha = 0$. That is,
$\alpha$ can be considered as an integrable 1-form in the space of variables 
$(x_p,...,x_n)$, parametrized by $(x_1,...,x_{p-1})$. We will distinguish two subcases: 
$d'\alpha = 0$ and $d'\alpha \neq 0$.

Subcase a). Suppose that $d'\alpha = 0$. Then according to Poincar\'e Lemma we have
$\alpha_j = \sum_{i=1}^{p-1} \alpha_j^i x_i + \partial / \partial x_j q^{(2)}$, where
$q^{(2)}$ is a quadratic function in the variables $(x_p,...,x_n)$. By a linear change of
coordinates on $(x_p,...,x_n)$, we have $q^{(2)} = \sum_{j=p}^{p+r} \pm x_j^2 / 2$,
for some number $r \geq -1$, and accordingly
$$ \alpha = \sum_{j=p}^{p+r}( \pm x_j + \sum_{i=1}^{p-1} \alpha_j^i x_i) dx_j +
\sum_{i=\overline{1,p-1}, j= \overline{p+r+1,n}} \alpha_j^i x_i dx_j $$
By a linear change of coordinates $(x_1,...,x_{p-1})$ on one hand, and
$(x_{p+r+1},...,x_{n})$ on the other hand, we can normalize the second part of the
above expression to obtain
$$\alpha = \sum_{j=p}^{p+r}( \pm x_j + \sum_{i=1}^{p-1} \tilde{\alpha}_j^i x_i ) +
\sum_{j=1}^s x_j dx_{p+r+j}$$
for some number $s$ ($0 \leq s \leq min(p-1, n-p-r)$).
                                                                            
Replacing $x_j$ ($j= \overline{p,p+r}$) by new $x_j = x_j \mp  \tilde{\alpha}_j^i x_i$ we have
$\omega = dx_1 \wedge ... \wedge dx_{p-1} \wedge \alpha$
where
$$ \alpha = d [ \sum_{j=p}^{p+r} \pm x_j^2 / 2 + \sum_{i=1}^s x_i x_{p+r+i}] $$
(with $- 1 \leq r \leq q = n-p, 0 \leq s \leq q - r$). These are the linear co-Nambu forms
of Type 1 in Theorem \ref{thm:Linear}.

Subcase b). Suppose that $d'\alpha \neq 0$. Then since $d'\alpha$ is of constant 
coefficients, we can change the coordinates $(x_p,...,x_n)$ linearly so that
$d'\alpha = dx_p \wedge dx_{p+1} + ... + dx_{p+2r} \wedge dx_{p+2r+1}$
in these new coordinates, for some $r\geq 0$.

If $r \geq 1$, then considering the coefficients of the term 
$dx_p \wedge dx_{p+1} \wedge dx_{i} \, (i > p+1), dx_p \wedge dx_{p+2} \wedge dx_{p+3}$
and $dx_{p+1} \wedge dx_{p+2} \wedge dx_{p+3}$ in $0 = \alpha \wedge d'\alpha$, we obtain
that all the coefficients of $\alpha$  are zero, i.e. $\alpha = 0$, which is absurd.
Thus $d'\alpha = dx_p \wedge dx_{p+1}$, and the condition  $\alpha \wedge d'\alpha = 0$
implies $\alpha = \alpha_1 dx_p + \alpha_2 dx_{p+1}$ with linear functions $\alpha_1$
and $\alpha_2$ depending only on $x_1,...,x_{p-1}, x_p, x_{p+1}$. In this Subcase b),
$\omega = dx_1 \wedge ... \wedge dx_{p-1} \wedge \alpha$ also has the form
(\ref{eq:TempoType2}), as in Case 2.

Suppose now that $\omega$ has the form (\ref{eq:TempoType2}), 
as in Case 2 or Subcase b) of Case 1:
$$
\omega = \sum x_i \omega_i = \sum_{j=1}^{p+1} a_j dx_1 \wedge ... \wedge dx_{j-1}
\wedge dx_{j+1} \wedge ... \wedge dx_{p+1}
$$
There are also 2 subcases:

a) $\partial a_j / \partial x_i = 0 \; \; {\rm for}\; j=\overline{1,p+1},i=\overline{p+2,n}$. In other words,
$$ \omega = \sum_{i,j=1}^{p+1} a^i_j x_idx_1 \wedge ... \wedge dx_{j-1}
\wedge dx_{j+1} \wedge ... \wedge x_{p+1} $$
with constant coefficients $a^i_j$. 

To see that $(a^i_j)$ can be put in Jordan form, notice that the linear Nambu tensor
corresponding to $\omega$ is, up to multiplication by a constant:
$$ \Pi = (\sum_{i,j=1}^{p+1} \pm a^i_j x_i \partial / \partial x_j) \wedge
\partial / \partial x_{p+2} \wedge ... \wedge \partial / \partial x_n $$
The first term in $\Pi$ is a linear vector field, which is uniquely defined by a linear
transformation ${\Bbb R}^{p+1} \to {\Bbb R}^{p+1}$ given by the matrix $(a^i_j)$, 
so this matrix can be put in Jordan form.

b) There is $j \leq p+1$ and $i \geq p+2$ such that $\partial a_j / \partial x_i \neq 0$.
We can assume that $\partial a_1 / \partial x_n \neq 0$. Putting 
$A = \partial/ \partial x_3 \wedge ... \wedge \partial / \partial x_{p+1}$ in
$0 = i_A \omega \wedge d\omega$ we obtain
$$
0 = (a_1 dx_2 + a_2 dx_1) \wedge \sum_{i=\overline{1,n}, j=\overline{1,p+1}} dx_i \wedge
\frac{\partial a_j}{\partial x_i} 
dx_1 \wedge ... \wedge dx_{j-1} \wedge dx_{j+1} \wedge ... \wedge dx_{p+1}
$$ 
Considering the coefficient of $dx_1 \wedge ... \wedge dx_{p+1} \wedge dx_n$
in the above equation, we have 
$$ a_1 \partial a_2 / \partial x_n - a_2 \partial a_1 / \partial x_n = 0$$
Since $\partial a_1 / \partial x_n \neq 0$, it follows that $a_2$ is linearly dependent of
$a_1$. Similarly, $a_j$ is linearly dependent of $a_1$ for any $j=\overline{1,p+1}$. Thus
$\omega = a_1 \omega_1$ where $\omega_1$ is decomposable and constant:
$\omega_1 = dx_1 \wedge ... \wedge dx_p$ in some linear system of coordinates.
If $a_1$ is linearly independent on $(x_1,...,x_p)$ then we can also assume that
$a_1 = x_{p+1}$. Thus also in this Subcase b), $\omega$ is of Type 2 in Theorem
\ref{thm:Linear}.
$\Box$

The form of $\omega$ gives us a clear picture about the singular foliations associated to
linear Nambu structures: The foliation of a linear Nambu structure of Type 1 has 
$p$ first integrals, namely  $x_1, ..., x_{p-1}$ and $\sum_{j=p}^{p+r}\pm x_j^2 +
\sum_{j=1}^s x_j x_{p+r+j}$, and the leaves of the foliation are uniquely determined
by these first integrals. The singular foliation of a linear Nambu structure of Type 2
is a Cartesean product of a foliation given by a linear vector field in a linear space with
(an 1-leaf foliation on) another linear space.

Rewriting Theorem \ref{thm:Linear} in terms of Nambu tensors, we have:

\begin{cor}
\label{cor:LinearTensor}
Every linear Nambu tesnsor $\Pi$ of order $q = n-p \geq 3$ on an
$n$-dimensional linear space $V$ belongs to one of the following two types:

\underline{Type 1}:
$
\Pi = \sum_{j=1}^{r+1} \pm x_j \partial / \partial x_1 \wedge ... \wedge 
\partial / \partial x_{j-1}	\wedge 
\partial / \partial x_{j+1}	\wedge ... \wedge
\partial / \partial x_{q+1} +    \\
\sum_{j=1}^{s} \pm x_{q+1+j} \partial / \partial x_{1} \wedge ... \wedge
\partial / \partial x_{r+j} \wedge 
\partial / \partial x_{r+j+2} \wedge 
\partial / \partial x_{q+1}
$
(with $-1 \leq r \leq q, 0 \leq s \leq min(p-1, q-r)$).

\underline{Type 2}:
$
\Pi = 
\partial / \partial x_{1} \wedge ... \wedge
\partial / \partial x_{q-1} \wedge
(\sum_{i,j=q}^{n} b^i_j x_i \partial / \partial x_j)
$
\end{cor}

{\it Remark}. Linear Nambu tensors are the same as finite-dimensional Nambu-Lie
algebras (cf. \cite{Gautheron,Takhtajan}). Thus Corollary \ref{cor:LinearTensor}
can be viewed as the classification of finite-dimensional Nambu-Lie algebras. 
The case of 4-dimensional Nambu-Lie algebras of order 3 has been done in 
\cite{Gautheron}.

We notice here a very interesting duality between Type 1 and Type 2: The formula for $\Pi$
of Type 1 looks similar to that for $\omega$ of Type 2, and vice versa. 
This duality will play an important role in the rest of this paper. We should
notice also that if a differential form $\omega$ can be written in one of the two
forms presented in Theorem \ref{thm:Linear}, then it is obviously a linear co-Nambu
form.

We have the following natural notion of nondegeneracy for linear Nambu structures:

\begin{define}
\label{define:Nondegenerate}
A linear co-Nambu $p$-form $\omega$ (and its corrsponding linear Nambu $q$-tensor $\Pi$) 
of Type 1 is called {\bf nondegenerate} if and only if it can be written in the
form $ \omega = dx_1 \wedge ... \wedge dx_{p-1} \wedge d q^{(2)}$, where 
$q^{(2)} = \sum_{p}^{n} \pm x_j^2$ (is nondegenerate). In this case, $\omega$ and
$\Pi$ are called {\bf elliptic} if $q^{(2)}$ is negative-definite
or positive-definite. The absolute value of the signature of the quadratic function 
$q^{(2)}$ is called the {\bf signature} of $\omega$. 
The {\bf index} of $\omega$ is the index of $q^{(2)}$, defined only up to the 
involution $m \mapsto q+1-m$.   

A linear co-Nambu $p$-form $\omega$ (and its corresponding linear Nambu $q$-tensor $\Pi$) 
of Type 2 is called {\bf nondegenerate} if and only if it can be written in the
form $ \omega = \sum_{i=1}^{p+1} \sum_{j=1}^{p+1} a_i^j x_j dx_1 \wedge ... 
\wedge dx_{i-1} \wedge dx_{i+1} \wedge ... \wedge dx_{p+1}$, with $(a^i_j)$ being
nondegenerate, i.e. having non-zero determinant.
\end{define}

It is evident that a linear Nambu structure of Type 1 is nondegenerate if and only
if all the other linear Nambu structures nearby it are equivalent
to it in a natural sense, and there is only a finite number of equivalence
classes in this case, which are classified by the signature of $q^{(2)}$. On the other
hand, for nondegenerate linear Nambu structures of Type 2, there is a continuum
of equivalence classes, which are classified by the Jordan form of $(a^i_j)$, 
modulo multiplication by a non-zero number.

\section{Decomposition of nondegenerate Nambu singularities}
\label{section:Decomposition}

We will say that a singularity of a Nambu structure is {\it of Type 1} 
({\it of Type 2, nondegenerate, elliptic, hyperbolic}) if its linear part is so.
In this Section we will show that  Nambu structures are 
decomposable also at nondegenerate singularities.

\begin{thm}
\label{thm:Decomposition}
a) Let $O \in V$ be a nondegenerate singular point of Type 1 of a 
co-Nambu p-form (of co-order $q \geq 3$) $\omega$. Then in a small neighborhood
of $O$ in $V$, $\omega$ is decomposable: it can be written as
$$ \omega = \gamma_1 \wedge ... \wedge \gamma_{p-1} \wedge \alpha $$
where $\gamma_i$ are 1-forms which do not vanish at $O$, and $\alpha$ is an
1-form which vanishes at $O$.

b) Let $O \in V$ be a nondegenerate singular point of Type 2 of a 
Nambu q-tensor (of order $q \geq 3$) $\Pi$. If $q=n-1$ then we will also assume
that in the normal form of its linear part
$
\Pi^{(1)} =  
\partial / \partial x_{1} \wedge ... \wedge                      
\partial / \partial x_{q-1}
\wedge
(\sum_{i,j=q}^{n} b^i_j x_i \partial / \partial x_j)
$ 
as given in Corollary \ref{cor:LinearTensor}, the ($2 \times 2$) matrix
$(b^i_j)$ has a non-zero trace. Then in a small neighborhood                               
of $O$ in $V$, $\Pi$ is decomposable: it can be written as 
$$ \Pi = V_1 \wedge ... \wedge V_{q-1} \wedge X $$
where $V_i$ are vector fields which do not vanish at $O$, and $X$
is a vector field which vanishes at $O$.
\end{thm}

{\it Proof}. First we will prove a). The proof will not make use of the integrability
of Nambu tensors (or similar property of co-Nambu forms), so in fact the above theorem
can be stated in a stronger form.

According to the definition of nondegenerate singularities of Type 1, we
can suppose that $\omega$ has a Taylor expansion $\omega = \omega^{(1)} + 
\omega^{(2)} + ...$, with $\omega^{(1)} = dx_1 \wedge ... \wedge dx_{p-1} \wedge dq^{(2)}$,
where $q^{(2)} = \sum_{j=p}^n \pm x_j^2 / 2$. Express $\omega$ as a polynomial in
$dx_1, ..., dx_{p-1}$:
$ 
\omega = dx_1 \wedge ... \wedge dx_{p-1} \wedge \alpha + 
\sum_{j=1}^{p-1} dx_1 \wedge ... \wedge dx_{j-1} \wedge dx_{j+1}
\wedge ... \wedge dx_{p-1} \wedge \beta_j + 
\sum_{1 \leq i < j \leq p-1} dx_1 \wedge ... \wedge dx_{i-1} \wedge dx_{i+1}
\wedge ... \wedge dx_{j-1} \wedge dx_{j+1}
\wedge ... \wedge dx_{p-1} \wedge \gamma_{ij} + ... 
$
Here $\alpha, \beta_i, \gamma_{ij},...$ are differential forms which, when written in
coordinates $(x_1,...,x_n)$, do not contain the terms $dx_1,..., dx_{p-1}$. Applying
the equation $i_A \omega \wedge \omega = 0$ to $A = \partial / \partial x_1 \wedge ...
\wedge \partial / \partial x_{p-1}$, we have $\alpha \wedge \omega = 0$. In follows that
$\alpha \wedge \beta_j = 0, \alpha \wedge \gamma_{ij} = 0$, etc. We can consider
$\alpha$ and $\beta_j$ as differential forms on the space of variables $\{x_p, ..., x_n\}$,
parametrized by $x_1,..., x_{p-1}$, and by our assumption of nondegeneracy, we can apply
DeRham division theorem (cf. \cite{DeRham}), which says that, since the number of variables
is $q + 1  > 2$ which is the order of $\beta_j$, $\beta_j$ is divisable by $\alpha$:
$\beta_j = \alpha \wedge \theta_j$ where $\theta_j$ are smooth 1-forms.

Applying the equation $i_A \omega \wedge \omega = 0$ to 
$A = \partial / \partial x_1 \wedge ... 
\wedge \partial / \partial x_{j-1} \wedge 
\partial / \partial x_{j+1} \wedge ...
\wedge \partial / \partial x_{p-1} \wedge \partial / \partial x_{p}$, we get
$$
0= \omega \wedge [<\alpha, \partial / \partial x_p> ( (-1)^{p-j} dx_j + \theta_j ) - 
<\theta_j, \partial / \partial x_p> \alpha]
$$
Since $<\alpha, \partial / \partial x_p> = <\alpha^{(1)}, \partial / \partial x_p> + ... =
\pm x_p + ... \neq 0$, and we already have $\omega \wedge \alpha = 0$, we get that
$ \omega \wedge \gamma_j = 0 $
where $\gamma_j = dx_j + (-1)^{p-j} \theta_j$
Since $\gamma_j$ do not vanish are are linearly independent at $O$, it follows that
$\omega$ is divisible by the product of $\gamma_j$:
$$ \omega = \gamma_1 \wedge ... \wedge \gamma_{p-1} \wedge \alpha' $$
for some 1-form $\alpha'$. By adding a combination of $\gamma_j$ to $\alpha'$,
we can assume that $\alpha'$ does not contain the terms $dx_1,...,dx_{p-1}$ when written
in the coordinates $(x_1,...,x_n)$. Then considering the terms containing
$dx_1 \wedge ... \wedge dx_{p-1}$ on the two sides of the equation
$ \omega = \gamma_1 \wedge ... \wedge \gamma_{p-1} \wedge \alpha' $, 
it follows that in fact we have $\alpha' = \alpha$. Statement a) is proved.

The proof of Statement b) in case $q \leq n-2$ is the same as that of a), 
by the duality $vector \leftrightarrow covector$. We will now prove b) for
the case $q = n-1$. In this case we have
$                                                                         
\Pi = 
\partial / \partial x_{1} \wedge ... \wedge                             
\partial / \partial x_{n-2}
\wedge (X_{n-1} \partial / \partial x_{n-1} + X_n \partial / \partial x_{n}) +
(\sum_{i=1}^{n-2} B_i  \partial / \partial x_{1} \wedge ... \wedge  
\partial / \partial x_{i-1} \wedge \partial / \partial x_{i+1} \wedge ... \wedge                    
\partial / \partial x_{n-2})  \wedge  
\partial / \partial x_{n-1} \wedge \partial / \partial x_{n}
$,                                                                         
where $B_i$ contains only terms of degree $ \geq 2$ in the Taylor expansion,
and the linear part of  the vector field
$X = X_{n-1} \partial / \partial x_{n-1} + X_n \partial / \partial x_{n}$ has
non-zero trace, that is 
$\partial X_{n-1 }/ \partial x_{n-1} + \partial X_n / \partial x_{n} \neq 0$.
Notice that $X$ is a Hamiltonian vector field of $\Pi$, given by 
the $(q-1)$-tuple of functions $(x_1,...,x_{n-2})$ (here $n-2 = q-1$). Hence
$X$ preserves $\Pi$: ${\cal L}_X \Pi = 0$. Considering the coefficient
of the term 
$
\partial / \partial x_{1} \wedge ... \wedge               
\partial / \partial x_{i-1} \wedge \partial / \partial x_{i+1} \wedge ... \wedge 
\partial / \partial x_{n-2}                                            
$
in the equation   ${\cal L}_X \Pi = 0$, we obtain a relation of the form
$$ 
(X(B_i) - 
(\partial X_{n-1 }/ \partial x_{n-1} + \partial X_n / \partial x_{n}) B_i)
\partial / \partial x_{n-1} \wedge \partial / \partial x_{n} 
= U \wedge X
$$
for some 
$U = U_{n-1} \partial / \partial x_{n-1} + U_n \partial / \partial x_{n}$
Since $\partial X_{n-1 }/ \partial x_{n-1} + \partial X_n / \partial x_{n} \neq 0$,
it follows that we have a relation of the form
$B_i = V_{n-1}X_n - V_n X_{n-1}$,
and therefore 
$B_i \partial / \partial x_{n-1} \wedge \partial / \partial x_{n}$
is divisible by $X$:
$
B_i \partial / \partial x_{n-1} \wedge \partial / \partial x_{n} =  
(V_{n-1} \partial / \partial x_{n-1} + V_n \partial / \partial x_{n}) \wedge  X 
$
Thus in this case, by using the fact that $X$ has non-zero trace, instead of
its nondegeneracy, we also obtain the divisibility by $X$ of the terms of degree
$q-2$ in the expression of $\Pi$ as a polynomial in 
$\partial/\partial x_1,... \partial / \partial x_{q-1}$. The rest of the 
proof is the same as for the case $q \leq n-2$.
$\Box$

The nondegeneracy implies that the 1-form $\alpha$ in the above theorem, considered
as an 1-form on the space of the variables $(x_p,...,x_n)$, will have exactly one
(nondegenerate) zero point for each value of the parameters $(x_1,...,x_{p-1})$, 
and of course this zero point will depend smoothly on the parameters $(x_1,...,x_{p-1})$.
A similer statement is true for the vector field $X$ in the second case. Thus we have:

\begin{cor}
\label{cor:ZeroSet}
If $O$ is a nondegenerate singular point of Type 1 of a Nambu tensor $\Pi$ of order
$q \geq 3$ in an $n$-dimensional manifold, then the set of zero points of $\Pi$ near 
$O$ forms a $(n-q-1)$-dimensional submanifold. 
If $O$ is a nondegenerate singular point of Type 2 of a Nambu tensor $\Pi$ of order
$q \geq 3$ in an $n$-dimensional manifold (when $q= n-1$ 
we need the same additional assumption as in the previous
theorem), then the set of zero points of $\Pi$ near $O$ forms a 
$(q-1)$-dimensional submanifold. 
\end{cor}

\section{Nondegenerate singularities of Type 1}
\label{section:Type1}

We have the following result about the linearization of co-Nambu forms of Type 1:

\begin{thm}
\label{thm:FormType1}
Let $O$ be a nondegenerate singular point of Type 1 of a smooth co-Nambu $p$-form 
$\omega$ of co-order $q > 2$.

a) If the singular point $O$ is of elliptic type then $\omega$ is linearizable 
in a neighborhood of $O$, up to multiplication by a non-zero smooth function. 
In other words, there is a local smooth system of coordinates $(x_1,...,x_n)$ 
in a neighborhood of $O$ 
such that we have $\omega =  f dx_1 \wedge ... \wedge dx_{p-1} \wedge 
\alpha^{(1)}$, where $f$ is a smooth function which does not vanish at $O$, 
and $\alpha^{(1)} = dq^{(2)}$ is a nondegenerate linear closed 1-form in the
variables $(x_p,...,x_n)$ (which does not depend on $(x_1,...,x_{p-1})$).

b) If $\omega$ is analytic (real or complex), then it is linearizable analytically
in a neighborhood of $O$, 
up to multiplication by an analytic function 
wich does not vanish at $O$.

c) If $\omega$ is only $C^{\infty}$ but not analytic, and $O$ is not of 
elliptic type, then $\omega$ is still formally linearizable at $O$, up to 
multiplication by a formal function which does not vanish at $O$.
\end{thm}

{\it Proof}. Statement a) and Statement b) of the above theorem are
absolutely similar to that of Reeb \cite{WR}, as improved by Moussu 
\cite{Moussu}, for the case of integrable 1-forms, and the proof is essentially
the same except for some additional regular first integrals. So we will only give a 
sketch of the proof here. The details of the steps can be found in 
\cite{Moussu,WR}. In the elliptic case, we can blow up along the local 
$(p-1)$-dimensional 
submanifold of elliptic singular points of $\omega$ (cf. Corollary 
\ref{cor:ZeroSet}), and then take a double covering of the blown-up manifold. 
In this double covering we have a regular foliation induced by the foliation 
associated to the Nambu structure. All the leaves of this foliation are 
diffeomorphic to $S^{q}$ due to Reeb's stability theorem, and the foliation 
itself is a regular fibration of fiber $S^{q}$. On the $p$-dimensional
base space of this fibration we have a smooth involution, whose fixed point set
is a local $(p-1)$-dimensional manifold (which corresponds to the manifold of 
zero points of $\omega$). It follows that there is a system of coordinates
$(f_1,...,f_p)$ on the base manifold of the fibration such that 
$(f_1,...,f_{p-1})$ are invariant under the involution, $f_p = 0$ on the 
submanifold of fixed points and $f_p^2$ is invariant under the involution. 
These coordinates give rise to the first integrals of the singular foliations 
of the Nambu structure: the first $(p-1)$ first integrals are regular and 
functionally independent, the last one is zero on the submanifold of zero points 
of the co-Nambu form $\omega$ and is nondegenerate positive-definite in the 
transversal direction to this submanifold. Taking the first $(p-1)$ first 
integrals as coordinates and applying the Morse's lemma
to the last first integral, we get the linearization of $\omega$ up to 
multiplication by a non-zero smooth function.
In the real analytic case, one can complexify the picture, then realify it 
back (in a different way) so that the singularity becomes elliptic, and then 
proceed as above. The complex analytic case is similar, without the step of 
complexifying.

Let us now prove Statement c) of the theorem. Since in this case the blowing-up
argument does not work so easily, we adopt a different strategy.
By induction we assume that we have found a new system of local
coordinates $(x_1,...,x_n)$ such that the Taylor expansion of $\omega$ in these
coordinates have ``good'' (r-1) first terms:
$$
\begin{array}{l}
\omega^{(1)} = dx_1 \wedge ... \wedge dx_{p-1} \wedge \alpha^{(1)} \\
\omega^{(2)} = dx_1 \wedge ... \wedge dx_{p-1} \wedge \alpha^{(2)} \\   
...... \\
\omega^{(r-1)} = dx_1 \wedge ... \wedge dx_{p-1} \wedge \alpha^{(r-1)} 
\end{array}
$$
where $\alpha^{(1)} = \sum_p^n \pm x_i dx_i$.
When $r=2$, this assumption follows from the definition of nondegenerate 
singularities of Type 2. We will show that we can make so that
$$\omega^{(r)} = dx_1 \wedge ... \wedge dx_{p-1} \wedge \alpha^{(r)}$$

Let us use the following notations:
\begin{equation}
\label{eq:Notations}
\begin{array}{l}
x = (x_1,...,x_{p-1}) \\
y = (x_p,..., x_n) = (y_1,..., y_{q+1}) \\
dx = dx_1 \wedge ... \wedge dx_{p-1} \\
d\hat{x_i} = dx_1 \wedge ... \wedge dx_{i-1} \wedge dx_{i+1} 
\wedge ... \wedge dx_{p-1} \\
\partial x =  \partial / \partial x_1 \wedge ... \wedge
\partial / \partial x_{p-1} \\
\partial \hat{x_i} = \partial / \partial x_1 \wedge ... \wedge
\partial / \partial x_{i-1} \wedge \partial / \partial x_{i+1} \wedge
\partial / \partial x_{p-1}
\end{array}
\end{equation}

Decompose $\omega^{(r)}$ into 
$\omega^{(r)} = dx \wedge \alpha^{(r)}+ \omega'$ , where
$\omega'$ consists of the terms which are not divisible by
$dx$. Put $A= \partial \hat{x_k} \wedge \partial / \partial y_1$  
for some index $k < q$. We have that
$i_A\omega^{(r)} = \sum_{j} \nu_j dy_j \pm \alpha^{(r)}_1 dx_k$ for some 
functions $\nu_j$. The terms of degree $r+1$ in the relation $i_A \omega \wedge
\omega = 0$ give:
\begin{equation}
\label{eq:*}
\pm y_1 dx_k \wedge \omega' + (\sum_{j} \nu_j dy_j) \wedge dx
\wedge \alpha^{(1)} = 0,
\end{equation}
which implies that
$
\pm y_1 dx_k \wedge \omega' = dx \wedge \gamma_k
$
for some $\gamma_k$, and
$
x_k \wedge \omega' = dx \wedge \gamma_k'
$
for some $\gamma_k'$.
By varying $k$ from 1 to $p-1$, we obtain that
\begin{equation}
\label{eq:**}
\omega' = \sum_{k=1}^{p-1} d\hat{x_k} \wedge \omega^{k}
\end{equation}
with $\omega^k = \sum_{i < j}\omega^k_{ij} dy_i \wedge dy_j$.

Putting Equation (\ref{eq:**}) into the left-hand side of Equation (\ref{eq:*})
we get
$$
\pm y_1 \omega^k \wedge dx                              
= (\sum_{j} \nu_j d y_j) \wedge \alpha^{(1)} \wedge dx,  
$$
which implies
\begin{equation}
\label{eq:***} 
\omega^k = \alpha^{(1)} \wedge \beta^k
\end{equation}
with $\beta^k = \sum_{l}\beta^k_l dy_l$ for some $\beta^k_l$.

The term of degree $r$ in $i_A \omega \wedge d\omega = 0$ gives
$$
i_A\omega^{(1)} \wedge d\omega^{(r)} + i_A\omega^{(2)} \wedge d\omega^{(r-1)} +
... + i_A\omega^{(r)} \wedge d\omega^{(1)} = 0,
$$
which implies $ \pm y_1 dx \wedge d\omega^k = 0$, and hence
$dx \wedge d \omega^k = 0$.

Thus the derivation of $\omega^k$ with respect to the variables $y$ is zero:
$d_y \omega^k = 0$. Putting the relation (\ref{eq:***}) in this equation we get
\begin{equation}
\label{eq:****}
\alpha^{(1)} \wedge d_y \beta^k = 0
\end{equation}

Now we will use the nondegeneracy of $\alpha^{(1)}$. The division theorem of
DeRham (cf. \cite{DeRham}) says that in this case we can divide $d_y \beta^k$
by $\alpha^{(1)}$: 
\begin{equation}
\label{eq:beta}
d_y \beta = \alpha^{(1)} \wedge \beta^{(r-2)}
\end{equation}
where $\beta^{(r-2)}$ is a homogeneous  1-form of degree $r-2$.
Differentiating (\ref{eq:beta}) with respect to the variables $y$, we get
$
\alpha^{(1)} \wedge d_y \beta^{(r-2)} = 0,
$
which implies
$$
 d_y \beta^{(r-2)} = \alpha^{(1)} \wedge \beta^{(r-4)}
$$
for some $\beta^{(r-4)}$.
Repeat the above process until we get a form $\beta^{(r-2h)}$ with 
$d_y \beta^{(r-2h)} = 0$. Then we go back: $\beta^{(r-2h)} = d_y \phi^{(r-2h+1)}$
and the equation $d_y \beta^{(r-2h+2)} =\alpha^{(1)} \wedge \beta^{(r-2h)}$ gives
$\beta^{(r-2h+2)} = -\phi^{(r-2h+1)} \alpha^{(1)} + d_y \phi^{(r-2h+3)}$. Keep 
going back until we refind $\beta^k$ in the form
$$ \beta^k = - \phi^{(r-1)}_k \alpha^{}(1) + d_y \phi^{(r+1)}_k$$
It follows that in (\ref{eq:***}) we can change $\beta^k$  by an exact 1-form:
$$ 
\omega^k = \alpha^{(1)} \wedge d_y \phi^{(r+1)}_k 
$$
Consider now the following new system of coordinates
\begin{equation}
\begin{array}{ll}
x_1' & = x_1 \pm \phi^{(r+1)}_1 \\
...  & ...                                 \\
x_{p-1}' & = x_{p-1} \pm \phi^{(r+1)}_{p-1}  \\
y' & = y
\end{array}
\end{equation}
In these new coordinates, $\omega^{(1)}$ becomes
$dx_1' \wedge ... \wedge dx_{p-1}' \wedge (\sum_{j=1}^{q+1} \pm y_j dy_j)
= \omega^{(1)} + \sum_{k=1}^{p-1} \pm d\hat{x_k} \wedge \alpha{(1)}
\wedge d_y \phi^{(r+1)_k} + \; (terms \; of \; degree \; > q)$

Thus, by choosing appropriate signs in the above change of variables, we can 
kill the term $\omega' = \sum d \hat{x_k} \wedge \omega^k = 
\sum_{k=1}^{p-1} \pm d \hat{x_k} \wedge \alpha{(1)} \wedge d_y \phi^{(r+1)}_k$ 
in the expression $\omega^{(r)} = dx \wedge \alpha^{(r)} + \omega'$

Repeating the above procedure for $r$ going from 2 to infinity, we find a 
formal system of coordinates $(x_1,...,x_n)$ in which $\omega = \sum_{r=1}^{\infty}
\omega^{(r)}$ with
$\omega^{(r)} = dx_1 \wedge ... \wedge dx_{p-1} \wedge \alpha^{(r)}$ 
for every $r$. In particular, 
$\omega = dx \wedge \alpha$,
where $dx = dx_1 \wedge ... \wedge dx_{p-1}$ and $\alpha = \sum \alpha^{(k)}$. 
The relation $i_{\partial x} \omega \wedge d \omega = 0$ implies that
$\alpha \wedge d_y \alpha = 0$, that is $\alpha$ can be considered as an 
integrable 1-form in the variables $y$, with a nondegenerate closed
linear part $\alpha^{(1)}$, and parametrized by $O$. It is well-known that
in this case $\alpha$ is formally linearizable up to 
multiplication by a formal function $f$ (see e.g. \cite{Moussu}).
Theorem \ref{thm:FormType1} is proved.
$\Box$
                                                                                   
{\it Remarks}. In the above theorem, we have linearization only up to               
multiplication by a function, because $\omega$ is not a closed form in general.
It is closed (outside singular points) only up to multiplication by a               
function. There is another simple proof of the analytic (and formal) case of 
the above theorem, which uses Theorem \ref{thm:Decomposition} and Malgrange's 
Frobenius theorem with singularities \cite{Malgrange}. 

The above theorem implies that a nondegenerate Nambu tensor of Type is
(maybe formally) linearizable up to multiplication by a function. In fact, at least 
formally, we can 
linearize it without the need of multiplication by a function:

\begin{thm}
\label{thm:LinearizationType1}
Let $O$ be a nondegenerate singular point of Type 1 of a smooth Nambu $q$-tensor 
$\Pi$, $q > 2$. Then $\Pi$ is formally linearizable at $O$: there is a formal 
system of coordinates $(x_1,...,x_n)$ such that
$$
\Pi \stackrel{formally}{=} \sum_{i=1}^{q+1} \pm x_i \partial / \partial x_1 
\wedge ... \wedge \partial / \partial x_{i-1} \wedge \partial / \partial x_{i+1}
\wedge ... \wedge \partial / \partial x_n
$$                                       
\end{thm}

{\it Proof}. According to Theorem \ref{thm:FormType1}, we can write 
$\Pi = f \Pi_1$ where 
$
\Pi_1 = \sum_{i=1}^{q+1} \pm x_i \partial / \partial x_1      
\wedge ... \wedge \partial / \partial x_{i-1} \wedge \partial / \partial x_{i+1} 
\wedge ... \wedge \partial / \partial x_n                                        
$

We want to change the coordinates $(x_1,...,x_{q+1})$ (and leave 
$x_{q+2},...,x_n$ unchanged) so that to make $f = 1$. 
We will forget about the parameters $(x_{q+2},...,x_n)$ and will assume for 
simplicity that $n=q+1$

Write $f= \sum f^{(r)}$ where $f^{(r)}$ is homogenous of order $r$ in 
$(x_1,...,x_{q+1})$. By a change of coordinates of the type
$x_1' = g x_1,..., x_{q+1}' = g x_{q+1}$, 
we can make $f^{(0)} = 1$. We assume now that we already 
have $f^{(1)} = ... = f^{(r-1)} = 0$ for some $r \geq 1$. We will show that 
there is a change of coordinates which changes $x_i$ by terms of degree $\geq 
r$, and which kills $f^{(r)}$. It amounts to find a vector field $X$ such that
$$
{\cal L}_X \Pi_1 = f^{(r)} \Pi_1
$$
where ${\cal L}$ denotes the Lie derivative. 
Consider the volume form 
$\Omega = dx_1 \wedge ... \wedge dx_{q+1}$. Then it is easy to see, by 
contracting $\Pi_1$ with $\Omega$, that the equation
${\cal L}_X \Pi_1 = f^{(r)} \Pi_1$ is equivalent to the equation
$$
dX(Q) = (f^{(r)} + div_{\Omega} X) dQ,
$$
where $Q = 1/2 \sum \epsilon_i x_i^2$ with $\epsilon_i = \pm 1$. In turn, this equation is 
equivalent to the following system of equations:
$$
\begin{array}{l}
div_{\Omega} X + f^{(r)} = d [2 Q F(Q)] / dQ \\
X(Q) = 2 Q F (Q)
\end{array}
$$
where $F$ is an unknown function.
Write $X = A + Y$, where $A = F(Q) \sum x_i \partial / \partial x_i$, and 
$Y$ is a vector field such that $Y(Q) = 0$. Then the above system of equations 
is equivalent to a system of the type
$$
\begin{array}{l}
Y(Q) = 0       \\
\beta(Q) + div_{\Omega} Y = f^{(r)}
\end{array}
$$
where $\beta$ is an unknown function. The equation $Y(Q) = 0$ is equivalent to 
the fact that $Y = \sum_{i < j} f_{ij} Y_{ij}$ where $Y_{ij} = \epsilon_i x_j \partial 
/ \partial x_i - \epsilon_j x_i \partial / \partial x_j$. For such an $Y$, we 
have $div_{\Omega} Y = \sum_{i < j} Y_{ij}(f_{ij})$. Denote by $J$ the set of 
homogenous polynomials of degree $r$. The solvability of the above system of 
equation follows from the following facts, which can be verified easily by 
choosing appropriate $f_{ij}$:

1. If a monomial $x^I = x_1^{I_1}...x_{q+1}^{I_{q+1}}$ has one of $I_i$ to be an 
odd number, then it belongs to $J$.

2. $Q^s$ is equivalent to $\lambda x_1^{2s}$ modulo $J$ for some non-zero number 
$\lambda$.

3. Any monomial $x^I = x_1^{I_1}...x_{q+1}^{I_{q+1}}$, with all $I_i$ even, is 
equivalent to $\lambda x_1^{\sum I_i}$ modulo $J$ for some number $\lambda$.

Thus the above system of equations can always be solved. The theorem is proved.
$\Box$

Suppose now that $\omega$ is of class $C^{\infty}$, is nondegenerate of Type 1 
at a zero point $O$, is not elliptic at $O$ 
but has an index different from 2 and $q-1$ (i.e. signature different from 
$q-3$, cf. Definition \ref{define:Nondegenerate}). Then the regular local leaves 
of the fibration associated to the linear part $\omega^{(1)}$ of $\omega$ are 
simply-connected (they are diffeomorphic to a direct product of a disk with 
a sphere of dimension different from 1). It follows from Reeb's stability 
theorem that the local regular leaves of $\omega$ are diffeomorphic to that of 
$\omega^{(1)}$. One can show easily in this case that the singular foliation 
associated to $\omega$ is homeomorphic to the one associated to the linear part
of $\omega$. (See \cite{Medeiros} for the case $p=1$). According to
Moussu \cite{Moussu,Moussu2}, if $\omega$ is of order 
$p=1$ (i.e. is an 1-form) and its index is different from 2 and $q-1$, or if 
its index is 2 but all of its leaves are closed except for a finite number of 
leaves which contain the origin in the limit, then 
it admits a smooth first integral, which means that $\omega$
is smoothly linearizable up to multiplication by a smooth 
function which does not vanish at $O$. 
We suspect that it is also true for the case $p > 1$.
If $\omega$ is of index 2 at 
$O$ and without the condition on the closedness of the leaves, 
then it may have no local first integral, (which implies in 
particular that it may not be linearizable up to multiplication by a function),
as the following example shows:

{\bf Example}. Consider the 2-form 
$
\omega = [dq + l(q) \alpha] \wedge [dx_3 + h(q) \alpha]
$
near $0$ in ${\Bbb R}^{3+k}$, where $q = x_1^2 + x_2^2 - y_1^2 - ... - y_k^2$, $\alpha = 
\frac{x_1 dx_2 - x_2 dx_1}{x_1^2 + x_2^2}$ 
is a singular closed 1-form, $l(q)$ and $h(q)$ are 
two flat functions in $q$ at $0$ such that $l(q) = h(q) = 0$ when $q \leq 0$.
The conditions on $l(q)$ and $h(q)$ assure that $\omega$ is a smooth 2-form 
whose linear part at 0 is $dq \wedge dx_4$. Near a point $P$ such that 
$x_1(P)^2 + x_2(P)^2 \neq 0$, we can write $\alpha = df$ fior some function $f$.
Thus near this point $\omega$ can be considered as a pull-back of a 2-form on a
3-dimensional space via the map $(x_1,x_2,x_3,y_1,...,y_k) \mapsto (q,f,x_3)$. Since 
any 2-form on a 3-dimensional space is a co-Nambu form and a pull-back of a 
co-Nambu form is also a co-Nambu form, it follows that $\omega$ is a co-Nambu 
form. When $h(q) = 0$ and $l(q) > 0$ for $q > 0$, the leaves of the singular foliation 
associated to $\omega$ will spiral toward the cones $(q=0, x_3 = constant)$. In
this case the folitation has only one local first integral (up to functional 
dependence), which is $x_3$. If $l(q) \geq 0$ and $h(q)$ is not identically 0 when
$q > 0$, then the leaves of the singular foliation associated to $\omega$ also 
drift in $x_3$, and if we chose $l,h$ well enough this phenomenon will prevent 
the foliation from having a non-trivial first integral. For example, we can make 
$l(q)$ and $h(q)$ vanish together at a series of points $q_i$ which tend to $0$. 
Near each point 
$q_i$ we make $h(q)$ vary from positive to negative an infinite number of 
times and chose $l(q)$ and $h(q)$ so that the drift in terms of
$x_3$ of a leaf passing via some point $x \in {\Bbb R}^4$ 
with $q_i < q(x) < q_{i-1}$ and spiraling 
inwards or outwards (i.e. the curve drawn by the value of $x_3$ 
of a point on this leave while this point is moving 
inwards or outwards), is contained in a small interval $[-\epsilon_i, + \epsilon_i]$ 
($\lim \epsilon_i = 0$) and spans this interval an infinite number of times. It
follows that this leave contains the leaves $(q= q_i, x_3 = \; constant \; \in 
[-\epsilon_i,+\epsilon_i])$ or the leaves $(q= q_i, x_3 = \; constant \; \in 
[-\epsilon_i,+\epsilon_i])$
in its limit. By invariance with respect to $\partial /
\partial x_3$ and $x_1 \partial / \partial x_2 - x_2 \partial / \partial x_1$ 
of our construction, any other leaf nearby this leaf will have the same 
property (with the interval $[-\epsilon_i,+\epsilon_i]$ replaced by an interval
$[-\epsilon_i + \delta ,+\epsilon_i + \delta]$). It follows that for any local continuous
function $f$ which is invariant on the leaves of the folitaion, there is an 
open set containing 0 in the boundary, in which $f$ is constant.

\section{Nondegenerate singularities of Type 2}
\label{section:Type2}

\begin{thm}
\label{thm:Prelinearization}
Let $O$ be a nondegenerate singular point of Type 2 of a                   
Nambu tensor $\Pi$ of order $q \geq 3$ on an $n$-manifold $V$, whose linear 
part has the form                                                           
$                                                   
\Pi^{(1)} =                                         
\partial / \partial x_{1} \wedge ... \wedge         
\partial / \partial x_{q-1} \wedge                  
(\sum_{i,j=q}^{n} b^i_j x_i \partial / \partial x_j)
$.      
If $q=n-1$ then we will also assume that the matrix $(b^i_j)$ has a non-zero 
trace. Then there is a local system of coordinates $(x_1,..., x_n)$ in which 
$\Pi$ can be written as 
$$
\Pi = f                                         
\partial / \partial x_{1} \wedge ... \wedge          
\partial / \partial x_{q-1} \wedge X                  
$$
where $f$ is a function with $f(O) \neq 0$ and 
$X = \sum_{i=q}^n c_i(x_q,...,x_n) \partial / \partial x_i$ is a vector field 
which does not depend on $(x_1,...,x_{q-1})$.
\end{thm}

{\it Proof}. Write                                                              
$\Pi = V_1 \wedge ... \wedge V_{q-1} \wedge X $                                 
as in Theorem \ref{thm:Decomposition}, where in some local system               
of coordinates $(x_1,...,x_n)$ we have                                          
$V_i = \partial / \partial x_i +                                                
\sum_{k=q,...,n} v_i^k \partial / \partial x_k$                                 
and $X$ does not contain terms in $\partial / \partial x_i, \; i <q$,           
and has  $\sum_{i,j=q}^{n} b^i_j x_i \partial / \partial x_j$                   
as its linear part. Using Corollary \ref{cor:ZeroSet}, we can, and will,        
make so that $X = 0$ on the submanifold $(x_q = ... = x_{n} = 0)$.              
                                                                                
The integrability of $\Pi$ implies that for any pair of indices                 
$i,j < q$, $[V_i,V_j] \wedge V_1 \wedge ... \wedge V_{q-1} \wedge X = 0$.                  
Considering the terms containing $\partial/ \partial x_1 \wedge ... \wedge      
\partial / \partial x_{q-1}$ in this equation, we get that                      
$[V_i, V_j] \wedge X = 0$. Notice that $[V_i,V_j]$ may be considered as a       
vector field on the space of the variables $(x_q,...,x_n)$, parametrized by the 
the parameters $(x_1,...,x_{q-1})$. The nondegeneracy of the linear part of $X$ 
allows us to use DeRham division theorem (cf. \cite{DeRham}), which says that   
$[V_i,V_j]$ is divisible by $X$: $[V_i,V_j] = g_{ij} X$ where $g_{ij}$ is some  
smooth function. Similarly, we have that $[V_i, X] = f_i X$ where $f_i$ is some 
smooth function. Using these properties, we will change $\Pi$, $V_i$ and $X$ so that   
$\Pi$ is only changed by multiplication by a non-zero function, 
the above relations still hold, but in addition $V_1$ commutes with              
$X, V_2,..., V_{q-1}$.                                                          

The equation $[V_1, X] = f_1 X$ implies that $[V_1, gX] = (V_1(g) + f_1 g) X$ 
for any function $g$. The equation $V_1(g) + f_1 g = 0$ can be solved locally because 
$V_1$ is non-zero at $O$. Replacing $X$ by $gX$ and $\Pi$ by $g\Pi$, we still 
have  $\Pi = V_1 \wedge ... \wedge V_{q-1} \wedge X $, but with $[V_1, X] = 0$.
Assume now that we already have $[V_1,X] = 0$. For $i > 1, i < q$ we have     
$[V_1, V_i + \gamma_i X] = (g_{1i} + V_1(\gamma_i)) X$. One can easily solve         
the equation $(g_{1i} + V_1(\gamma_i)) = 0$ to find a $\gamma_i$ such that           
$[V_1, V_i + \gamma_i X] = 0$. Replacing $V_i$ by $V_i + \gamma_i X$, we get         
that $V_1$ commutes with $V_2,...,V_{q-1},X$.                                        
                                                                                     
Assume now that we already have that $V_1$ commutes with  $V_2,...,V_{q-1},X$.       
In other words, everything is invariant with respect to $V_1$.                      
Make the same process as above but with $V_2$, in a way which is invariant           
with respect to $V_1$, we get that                                                   
$V_2,..., V_{q-1}, X$ can be changed so that $\Pi$ remains the same but              
$V_2$ becomes commuting with $V_3,...,V_{q-1},X$. Repeating the above process        
with $V_3, V_4, ...$. In the end we get a new family of vector fields $V_i$ and $X$  
whose product is $\Pi$ and which commute pairwise.                                   
                                                                                     
Since $V_i$ commute pairwise and are linearly independent, there is a new local      
system of coordinates $(x_1,...,x_n)$ such that in these coordinates we have         
$V_i = \partial/ \partial x_i$ for $i=1,...,q-1$. The fact that $X$                  
commutes with $V_i$ means that the coefficients of $X$ in these coordinates          
will not depend on $(x_1,...,x_{q-1})$. Of course, we can also assume that             
$X$ does not contain the terms $\partial/ \partial x_i, \;i=1,...,q-1$, since        
substracting these terms from $X$ will not change $\Pi$. Thus $X$ can be             
considered as vector field on the space of the variables $(x_q,...,x_n)$, which      
vanishes at the origin (and which does not depend on the parameters                  
$(x_1,...,x_{q-1})$). 
$\Box$

\begin{thm}
\label{thm:LinearizationType2}
Let $O$ be a nondegenerate singular point of Type 2 of a
Nambu tensor $\Pi$ of order $q \geq 3$ on an $n$-manifold $V$, whose linear
part has the form
$                                                       
\Pi^{(1)} =                                                   
\partial / \partial x_{1} \wedge ... \wedge             
\partial / \partial x_{q-1} \wedge                      
(\sum_{i,j=q}^{n} b^i_j x_i \partial / \partial x_j)    
$.                                                       
If the matrix $(b^i_j)$ is non-resonant, i.e. if its eigenvalues 
$(\lambda_1,...,\lambda_{p+1})$
do not satisfy any relation of the form 
$\lambda_i = \sum_{j=1}^{p+1} m_j \lambda_j$ with $m_j$ being non-negative integers
and $\sum m_i \geq 2$, then $\Pi$ is 
smoothly linearizable, i.e. there is a local smooth system of coordinates 
$(x_1,...,x_n)$ in a neighborhood of $O$, in which $\Pi$ 
coincides with its linear part:
$$
\Pi =
\partial / \partial x_{1} \wedge ... \wedge           
\partial / \partial x_{q-1} \wedge                    
(\sum_{i,j=q}^{n} b^i_j x_i \partial / \partial x_j)  
$$
The above linearization can be made analytic if $\Pi$ is analytic and the 
eigenvalues $\lambda_1,...,\lambda_{p+1}$ of $(b^i_j)$ 
satisfy the Bryuno's incommensurability condition:
there exist positive constants $C, \epsilon$ such that for any $(p+1)$-tuple of 
non-negative integers $(m_1,...,m_{p+1})$ with $\sum m_i \geq 2$
and any index $k \leq p+1$ we have
$|(\sum \lambda_i m_i) - \lambda_k| > C \exp (- (\sum m_i)^{1-\epsilon})$.
\end{thm}

{\it Proof}.
Using Theorem \ref{thm:Prelinearization}, we can write
$\Pi = V_1 \wedge ... \wedge V_{q-1} \wedge Y$, where 
$ Y = f X = f \sum_{i=q}^n c_i(x_q,...,x_n) \partial / \partial x_i$. 
(We will forget about the fact that $V_i = \partial / \partial x_i$).
If the linear part of $X$ satisfies the nonresonance condition then we 
can apply Sternberg's theorem \cite{Sternberg} to linearize $X$ smoothly, and if
it satisfies the Bryuno's incommensurability condition then we can apply 
Bryuno's theorem (see e.g. \cite{Bryuno,Martinet}) to linearize $X$ 
analytically in the analytic case. Thus in both case we can assume that $X$ 
is already linearized and normalized:
$X = \sum_{i=p}^n \lambda_{i-p+1} x_i \partial / \partial x_i$
We want now to change $V_i$ and $Y$ so that they become commuting and the 
relation $\Pi = V_1 \wedge ... \wedge V_{q-1} \wedge Y$ still hold (without 
multiplying $\Pi$ by a non-zero function). 

Similarly to the proof of Theorem \ref{thm:Prelinearization}, 
we have $[V_1, Y] = [V_1, f X] = V_1(f) X = f_1 Y$ with $f_1 = V_1(f)/f$.  
The equation $V_1 (g) + f_1 g = 0$ can be solved locally because $V_1 (O) \neq 
0$. This time we will solve it on the submanifold $(x_q = ... = x_n = 0)$ of 
zero points of $\Pi$. So let $g$ be a non-zero function which does not depend 
on $(x_q,...,x_n)$ and which satisfies $V_1 (g) + f_1 g = 0$ on the submanifold
$(x_q = ... = x_n = 0)$. Then $[V_1/g, gY] = h/g Y$ where $h$ is a function 
which vanishes on the submanifold $(x_q = ... = x_n = 0)$. Under the 
nonresonance condition, a theorem of Roussarie \cite{Roussarie} says that the equation 
$Y(\gamma) = h$, or equivalently, $X(\gamma) = h/f$, has a smooth solution $h$.
(Notice here an important fact that $X$ does not depend on $(x_q,...,x_n)$, 
which allows us to use Roussarie's theorem). In the analytic case, the equation
$X(\gamma) = (\sum_{i=1}^{p+1} \lambda_i x_{q-1+i}
\partial /\partial x_{q-1+i}) (\gamma) = h/f 
= \sum_{s_q,...,s_n} (h/f)_{s_q,...,s_n}(x_1,...,x_{q-1}) 
x_q^{s_q}...x_n^{s_n}$ has the formal solution  
$
\gamma = \sum_{s_q,...,s_n} \frac{1}{\sum_{i=1}^{p+1} \gamma_i s_{q-1+i}} 
(h/f)_{s_q,...,s_n} x_q^{s_q}...x_n^{s_n},
$                                                
which can be verified easily to converge near $O$, under the incommensurability
condition of Bryuno. 

With a smooth or analytic function $\gamma$ such that $Y(\gamma) = h$, we have
$[V_1/g + \gamma Y, g Y] = 0$. Thus we can change $V_1$ by $V_1/g + \gamma Y$ 
and $Y$ by $gY$ to obtain $[V_1,Y] = 0$. Of course, this change does not affect
$\Pi$. After that, we can change $V_2,...,V_{q-1}$ so that they commute with 
$V_1$, in the same way as in the proof of Theorem \ref{thm:Prelinearization}.

Thus we can make $V_1$ commute with $V_2,...,V_{q-1},Y$, without affecting 
$\Pi$. Just as in the proof of Theorem \ref{thm:Prelinearization}, by induction
we can make $V_1,...,V_{q-1},Y$ commute pairwise. Then we can put $V_i = 
\partial / \partial x_i$ in some new local system of coordinates, and can assume that
$Y$ does not contain the terms 
$\partial / \partial x_1,...,\partial / \partial x_{q-1}$. Then we can linearize
$Y$, using Bryuno's or Sternberg's theorem, to finish the linearization of $\Pi$.
$\Box$

{\it Remark}. The Bryuno's incommensurability condition in the above theorem 
can indeed be replaced by a weaker so-called $(\Omega)$-condition plus the 
nonresonance condition (see e.g. \cite{Bryuno,Martinet} for the $(\Omega)$-condition).

Talking about co-Nambu forms of Type 2, the above theorems show that such 
co-Nambu forms can be written locally as
$
\omega = f \omega_1
$
where $f$ is some non-zero function and $\omega_1$ is a $p$-form which
does not contain the terms $dx_1,...,dx_{q-1}$ and does not 
depend on the variables $x_1,...,x_{q-1}$, in some local system of coordinates 
$(x_1,...,x_n)$. In other words, $\omega_1$ is a pull-back of a $p$-form on 
a $(p+1)$-dimensional space under a projection ${\Bbb R}^n \to {\Bbb R}^{p+1}$.
Furthermore, $\omega_1$ can be made linear if $\omega$ safisfies some
nonresonance or incommensurability condition. If $d \omega (O) \neq 0$, then 
a result of Medeiros \cite{Medeiros} (called fundamental lemma for integrable 
p-forms) says that 
$\omega$ itself is the pull-back of a $p$-form on 
a $(p+1)$-dimensional space under a projection ${\Bbb R}^n \to {\Bbb R}^{p+1}$.
Let us give a proof of this fact, which is a slight simplification of the one 
given in \cite{Medeiros}: 

First of all, notice that if $\omega$ is a co-Nambu 
$p$-form, then $d \omega$ is a co-Nambu $(p+1)$-form. Indeed, the condition 
(\ref{eq:Integrable}) in Definition \ref{define:coNambu} is trivial for 
$d \omega$, and the condition (\ref{eq:Decomposable}) about the decomposability
is easily verified: near 
a non-zero point of $\omega$ we can write $\omega = f dx_1 \wedge ... \wedge 
dx_p$, which implies $d \omega = df \wedge dx_1 \wedge ... \wedge dx_p$. 
If $d \omega (0) \neq 0$ 
then a Nambu tensor dual to it is regular at $O$ and gives rise to a local 
regular foliation, denoted by ${\cal F}$. The tangent spaces of ${\cal F}$ are 
nothing but the spaces of vectors whose contraction with $d \omega$ is zero. 
Therefore if $Z$ is a vector tangent to ${\cal F}$ at a point $x$ near $O$
we have $i_Z d \omega(x) = 0$. If $\omega (x) \neq 0$ then we also obtain that
$i_Z \omega (x) = 0$, by using again
the presentation   $\omega = f dx_1 \wedge ... \wedge dx_p$. Since the set of 
non-zero points of $\omega$ is dense near $O$ (because $d \omega (0) \neq 0$), 
by continuity we get that for any
$Z$ tangent to ${\cal F}$, $i_Z \omega (x) = 0$ and  $i_Z d \omega(x) = 0$.
It means that $\omega$ is locally a 
pull-back of a form on the local base space of ${\cal F}$. QED.

\end{document}